\documentclass[oneside,english]{amsart}
\usepackage[utf8]{inputenc}
\usepackage{color}
\usepackage{babel}
\usepackage{url}
\usepackage{enumitem}
\usepackage{amsthm}
\usepackage{amssymb}
\usepackage{graphicx}
\usepackage[bookmarks=true,bookmarksnumbered=true,bookmarksopen=false,
 breaklinks=true,pdfborder={0 0 0},pdfborderstyle={},backref=false,colorlinks=true]
 {hyperref}
\hypersetup{pdftitle={Structural Stability in Piecewise Mobius Transformations},
 pdfauthor={Renato Leriche and Guillermo Sienra},
 pdfsubject={Holomorphic Dynamical Systems},
 pdfkeywords={Dynamical Systems, Holomorphic Dynamics, Picewise Transformations},
 linkcolor=mygreen,citecolor=red,urlcolor=blue,anchorcolor=blue,menucolor=gray,runcolor=yellow}

\makeatletter
\newcommand{\lyxrightaddress}[1]{
	\par {\raggedleft \begin{tabular}{l}\ignorespaces
	#1
	\end{tabular}
	\vspace{1.4em}
	\par}
}
\theoremstyle{definition}
    \ifx\thechapter\undefined
      \newtheorem{defn}{\protect\definitionname}
    \else
      \newtheorem{defn}{\protect\definitionname}[chapter]
    \fi
\theoremstyle{remark}
    \ifx\thechapter\undefined
      \newtheorem{rem}{\protect\remarkname}
    \else
      \newtheorem{rem}{\protect\remarkname}[chapter]
    \fi
\theoremstyle{plain}
    \ifx\thechapter\undefined
	    \newtheorem{thm}{\protect\theoremname}
	  \else
      \newtheorem{thm}{\protect\theoremname}[chapter]
    \fi
\theoremstyle{definition}
    \ifx\thechapter\undefined
      \newtheorem{example}{\protect\examplename}
    \else
      \newtheorem{example}{\protect\examplename}[chapter]
    \fi
\theoremstyle{plain}
    \ifx\thechapter\undefined
  \newtheorem{cor}{\protect\corollaryname}
\else
      \newtheorem{cor}{\protect\corollaryname}[chapter]
    \fi
\theoremstyle{plain}
    \ifx\thechapter\undefined
      \newtheorem{conjecture}{\protect\conjecturename}
    \else
      \newtheorem{conjecture}{\protect\conjecturename}[chapter]
    \fi

\@ifundefined{date}{}{\date{}}
\usepackage{mathabx}
\usepackage{graphicx}
\usepackage{amsaddr}

\definecolor{mygreen}{rgb}{0.0,0.6,0.2}

\newcommand{\N}{\mathbb N}

\newcommand{\R}{\mathbb R}
\newcommand{\C}{\mathbb C}
\newcommand{\Chat}{\widehat{\C}}

\newcommand{\D}{\mathbb{D}}

\newcommand{\Bc}{\mathcal B}

\newcommand{\Nc}{\mathcal N}
\newcommand{\Pc}{\mathcal P}
\newcommand{\Rc}{\mathcal R}

\newcommand{\Fix}{\mathrm{Fix}}
\newcommand{\Per}{\mathrm{Per}}
\newcommand{\PerA}{\mathrm{Per_{attr}}}
\newcommand{\PerR}{\mathrm{Per_{rep}}}
\newcommand{\PerN}{\mathrm{Per_{neu}}}
\newcommand{\PerP}{\mathrm{Per_{par}}}
\newcommand{\PerE}{\mathrm{Per_{ell}}}
\newcommand{\PerID}{\mathrm{Per_{id}}}
\newcommand{\PerG}{\mathrm{Per_{ghost}}}

\newcommand{\selfarrow}{\ensuremath{\rotatebox[origin=c]{90}{\ensuremath{\circlearrowleft}}}}

\AtBeginDocument{

}

\makeatother

\providecommand{\conjecturename}{Conjecture}
\providecommand{\corollaryname}{Corollary}
\providecommand{\definitionname}{Definition}
\providecommand{\examplename}{Example}
\providecommand{\remarkname}{Remark}
\providecommand{\theoremname}{Theorem}

\begin{document}
\title{Structural Stability in Piecewise Möbius Transformations}
\author{Renato Leriche$^{1,a,*}$ \& Guillermo Sienra$^{1,b}$}
\date{December 2024.}
\begin{abstract}
Structural stability of piecewise Möbius transformations (PMTs) is
examined from various perspectives. A result concerning structural
stability, restricted to the space of PMTs, is derived using hyperbolic
characteristics of the component functions and the pre-singularities
set, which facilitates a holomorphic motion. The analogous concept
of J-stability for rational maps is defined and analyzed for PMTs,
revealing some connections to general structural stability. The definitions
of hyperbolic and expansive PMTs are introduced, demonstrating that
they are not equivalent and that neither implies structural stability.
By synthesizing the previous results and analyses, sufficient conditions
for structural stability are established. Lastly, an example of structural
stability within the tent maps family, extended to the complex plane,
is presented.
\end{abstract}

\maketitle

\lyxrightaddress{$^{1}$Affiliation: Departamento de Matemáticas, Facultad de Ciencias,\\
Universidad Nacional Autónoma de México (UNAM), México.\\
$^{a}$E-mail: \href{http://r_lerichev@ciencias.unam.mx}{r\_lerichev@ciencias.unam.mx}\\
$^{b}$E-mail: \href{http://guillermosienra@ciencias.unam.mx}{guillermosienra@ciencias.unam.mx}\\
$^{*}$Corresponding author.}

\section*{Introduction}


A piecewise function in a space is defined by transformations restricted
to the components that belong to a finite partition of the space.
The study of the dynamics of piecewise maps arises from various contexts,
such as the interval exchange transformations (see for instance \cite{BreEtAl2010,GutEtAl2008,Via2006}),
the piecewise plane isometries (see \cite{AshFu2002,AshGoe2005,AshGoe2009,BosGoe2003,Buz2001,Dea2002,Dea2006,Goe1998,Goe1999,Goe2000,Goe2001,Goe2003,GoeQua2009})
and the piecewise contractions on $\R^{n}$ (see \cite{BruDea2009,CatEtAl2015}),
in addition to having applications in engineering and relations with
other areas of mathematics (see \cite{diBEtAl2007,Goe2003}).

The focus of this research work is the dynamics of \emph{piecewise
Möbius transformations} (abbreviated by its acronym as PMTs) in the
Riemann sphere (see \cite{Cru2005,LerSie2019} for other publications
about these maps). The most exciting connection between PMTs and other
areas of mathematics is that they emerge as the monodromy maps of
complex polynomial vector fields. These complex vector fields provide
a means of approaching Hilbert's problem 16 (still open), which deals
with the number and localization of limit cycles of real polynomial
vector fields (see \cite{Cru2005}). This connection is not addressed
in this paper, but it is expected that the results presented here
will be helpful for research on that complex vector fields.

A first study about stability and structural stability for PMTs is
worked on in \cite{LerSie2019}. In that paper, the associated group
generated by the component functions plays a central role. First,
if the limit set of the group does not intersect the boundary of the
domain partition and the component functions are fixed, continuous
deformations of the boundary induce continuous deformations of the
pre-singularities set as a compact set with the Hausdorff metric.
This continuity represents a form of stability; however, the structural
stability of the PMTs dynamics is not assured.

A second result in \cite{LerSie2019} indicates that if the boundary
of the partition is fixed, the associated group is structurally stable,
and the boundary of the partition is contained in a fundamental region
of the group, then the corresponding PMT is structurally stable in
the space of conformal automorphisms on the Riemann sphere.

In this paper, we will present sufficient conditions for the structural
stability of PMTs that are independent of the structural stability
of the associated group. To establish these conditions, we will define
and analyze hyperbolicity, $\alpha$-expansivity, and the analogous
concept of J-stability of rational functions in the Riemann sphere
for PMTs.

\section{Piecewise Möbius Transformations}

First of all, let us establish the basic definitions.
\begin{defn}
A \emph{piecewise Möbius transformation} (abbr.\emph{ PMT}) is a pair
$(P,F)$ where
\begin{itemize}
\item $P=\left\{ R_{k}\subset\Chat\right\} _{k=1}^{K}$ is a set of \emph{regions}
such that:
\begin{itemize}
\item Each $R_{k}$ is a non-empty open and connected set.
\item Each $\partial R_{k}$ is the union of piecewise smooth simple closed
curves.
\item $R_{k}\cap R_{j}=\emptyset$ if $k\neq j$.
\item $\bigcup_{k=1}^{K}\overline{R_{k}}=\Chat$.
\end{itemize}
\item $F:\Chat\,\selfarrow$, where each \emph{component function} $F|_{R_{k}}=f_{k}$
is the restriction of a conformal automorphism of $\Chat$ and $F$
is undefined in $\bigcup_{k=1}^{K}\partial R_{k}$.
\item $P$ is minimal in relation to $F$, that is, if $\overline{R_{k}}\cap\overline{R_{j}}\neq\emptyset$
and it is a union of curves, then $f_{k}\neq f_{j}$.
\end{itemize}
\end{defn}
\begin{rem}
$F:\Chat\,\selfarrow$ is a shorthand notation for $F:\Chat\rightarrow\Chat$.
\end{rem}
\begin{defn}
The\emph{ region of conformality} of a PMT $(\left\{ R_{k}\right\} _{k=1}^{K},\,F)$
is 
\[
R(F)=\bigcup_{k=1}^{K}R_{k}.
\]
\end{defn}
\begin{defn}
The \emph{discontinuity set} of a PMT $(\left\{ R_{k}\right\} _{k=1}^{K},\,F)$
is 
\[
B(F)=\partial R(F)=\bigcup_{k=1}^{K}\partial R_{k}.
\]
\end{defn}
\begin{rem}
Notice that the set $B(F)$ can be interpreted as the set of singularities
of $F$, since $F$ is not defined in such a set.
\end{rem}
A central construction to understand the dynamics of PMTs is the pre-singularities
set, as is it for meromorphic functions.
\begin{defn}
The \emph{pre-discontinuity set} of a PMT $F$ is
\[
\Bc(F)=\overline{\bigcup_{n\geq0}F^{-n}(B(F))}
\]

\begin{rem}
$\Bc(F)$ is the set of points that eventually lands in $B(F)$ under
$F$, or accumulation of those points. Then, if $z\in\Bc(F)$, there
exists $N\in\N$ such that $F^{N}(z)$ is undefined, or is an accumulation
point of such pre-singularities.
\end{rem}
\end{defn}
\begin{rem}
The set $\Bc(F)$ is alternatively called \emph{spiderweb} of $F$
and denoted $Spid(F)$ (see \cite{Cru2005}), because of its resemblance
to the spider's constructions in some cases. The analogous of this
set is called the \emph{exceptional set} or simply \emph{discontinuity
set} in the theory of bi-dimensional piecewise isometries (see \cite{Goe1996,Goe2003}).
\end{rem}
Analogously, as in holomorphic dynamics, it can be defined the set
with regular dynamics from the pre-singularities set.
\begin{defn}
The \emph{regular set} of a PMT $F$ is 
\[
\Rc(F)=\Chat-\Bc(F).
\]
\end{defn}
Another important set in the study of the dynamics of PMTs is the
pre-singularities accumulation set, called the $\alpha$-limit set.
\begin{defn}
The \emph{$\alpha$-limit set} of a PMT $F$ is 
\[
\alpha(F)=\Bc(F)-\bigcup_{n\geq0}F^{-n}(B(F)).
\]
\end{defn}
Also, it can be defined the $\omega$-limit set.
\begin{defn}
The \emph{$\omega$-limit set} of a PMT $F$ is $\omega(F)=\bigcup_{z\in\Rc(F)}\omega(z,F)$,
where $\omega(z,F)$ is the \emph{$\omega$-limit set} of $z$ under
$F$.
\end{defn}
\begin{rem}
The $\omega$-limit set is not always forward invariant nor is it
always backward invariant, since can occur $\omega(F)\cap(\mathcal{B}(F)-\alpha(F))\neq\emptyset$
as we will see later.
\end{rem}
Several results about the dynamics of PMTs have been obtained, they
can be thought as an extension of the dictionary of Sullivan (see
\cite{LerSie2019}). Below we state some of those results.

In what follows, let $F$ be a PMT.
\begin{thm}
(See \cite{Cru2005} and \cite{LerSie2019}.) $\Rc(F)$ is the set
where the family $\left\{ F^{n}\right\} _{n\in\N}$ is normal, and
$\Bc(F)$ is the set where the family $\left\{ F^{n}\right\} _{n\in\N}$
is not normal.
\end{thm}
\begin{thm}
$\Bc(F)$ is backward invariant, $\Rc(F)$ is forward invariant, and
$\alpha(F)$ is strictly backward invariant and forward invariant.
\end{thm}
\begin{proof}
We will only prove the assertions regarding $\alpha(F)$. For the
assertions about $\Bc(F)$ and $\Rc(F)$, see \cite{Cru2005} and
\cite{LerSie2019}.

Let $z\in\alpha(F)$.
\begin{enumerate}
\item Assume that $F^{-1}(z)\neq\emptyset$ and $F^{-1}(z)\nsubseteq\alpha(F)$.
It follows that $F^{-1}(z)\cap B(F)=\emptyset$ for all $z$, since
$F$ is undefined in $B(F)$. If $F^{-1}(z)\subset\Bc(F)-\alpha(F)$,
then $z\in\Bc(F)-\alpha(F)$, which leads to a contradiction. If $F^{-1}(z)\cap\Rc(F)\neq\emptyset$,
then $\left\{ F^{n}\right\} _{n\geq0}$ is normal at some $z_{0}\in F^{-1}(z)$
and at $z$, resulting in a contradiction. Thus, we conclude that
$F^{-1}(\alpha(F))\subset\alpha(F).$
\item Assume that $F(z)\notin\alpha(F)$. If $F(z)\in\Bc(F)-\alpha(F)$,
then $z\in\Bc(F)-\alpha(F)$, which is a contradiction. If $F(z)\in\Rc(F)$,
then $\left\{ F^{n}\right\} _{n\geq0}$ is not normal at $F(z)$ because
it is also not normal at $z$, leading to a contradiction. Therefore,
we have $F(\alpha(F))\subset\alpha(F).$
\item It is possible that $F^{-1}(z)=\emptyset$, which would imply $F(\alpha(F))\subsetneq\alpha(F)$.
However, $\alpha(F)\subset F^{-1}(\alpha(F))$ is always true by definition,
so using incise (1), we find that $F^{-1}(\alpha(F))=\alpha(F).$
\end{enumerate}
\end{proof}
\begin{thm}
$\overset{\circ}{\alpha(F)}=\emptyset$, where $\overset{\circ}{\alpha(F)}$
denotes the interior of $\alpha(F)$.
\end{thm}
\begin{proof}
Suppose $\overset{\circ}{\alpha(F)}\neq\emptyset$. Then there exists
an open set $U$ such that $U\subset\overset{\circ}{\alpha(F)}$.
Let $z\in U$; then there exists $N\geq0$ such that $F^{-N}(B)\cap U\neq\emptyset$.
Therefore, $B\cap F^{N}(U)\neq\emptyset$, a contradiction since $\alpha(F)$
is forward invariant and $\alpha(F)\cap B=\emptyset$ by definition.
\end{proof}
Since periodic points of PMTs are fixed points of Möbius transformations,
they can be categorized into \emph{attracting} (grouped in the set
$\PerA(F)$), \emph{repelling} ($\PerR(F)$), \emph{elliptic} ($\PerE(F)$)
and \emph{parabolic} ($\PerP(F)$). Additionally, there are periodic
points $z$ of period $n$ of a PMT $F$ for which there exists a
neighborhood $U$ of $z$ such that $f^{n}|_{U}$ is the identity
map in $U$ (grouped in $\PerID(F)$, and called periodic points \emph{of
identity}). Naturally, the set of \emph{neutral} or \emph{indifferent}
periodic points is $\PerN(F)=\PerE(F)\cup\PerP(F)\cup\PerID(F)$.
\begin{thm}
~
\[
\PerR(F)\cup\PerP(F)\subset\alpha(F)\subset\Bc(F),
\]
\[
\PerA(F)\cup\PerP(F)\cup\PerE(F)\cup\PerID(F)\subset\omega(F),
\]
 and 
\[
\PerA(F)\cup\PerE(F)\cup\PerID(F)\subset\Rc(F).
\]
\end{thm}
\begin{proof}
The family $\left\{ F^{n}\right\} _{n\geq0}$ is not normal in repelling
and parabolic periodic points; thus $\PerR(F)\cup\PerP(F)\subset\Bc(F)$.
However, since $\left\{ F^{n}\right\} _{n\geq0}$ is not completely
defined in $\Bc(F)-\alpha(F)$, it follows that $\PerR(F)\cup\PerP(F)\subset\alpha(F)$.

On the other hand, the family $\left\{ F^{n}\right\} _{n\geq0}$ is
normal in attracting, elliptic and identity periodic points; therefore
$\PerA(F)\cup\PerE(F)\cup\PerID(F)\subset\Rc(F)$. 

Given that the periodic points reside in their own $\omega$-limit
set, we have $\PerA(F)\cup\PerE(F)\cup\PerID(F)\subset\omega(F)$.

Finally, for parabolic periodic points $z$, there exists $w\in\Rc(F)$
such that $z\in\omega(w,F)$; hence $\PerP(F)\subset\omega(F)$.
\end{proof}
Since PMTs have a set of singularities, there are regular components
that can also exhibit an analogous behavior to the Baker domains of
meromorphic functions.
\begin{defn}
A point $z_{0}$ is a \emph{ghost-periodic} of \emph{period} $n$
of $F$ if $z_{0}\in F^{-N}(B)$ for some $N\geq0$ and there exists
a periodic regular component $U$ of period $n$ such that $z_{0}\in\partial U$
and for all $z\in U$
\[
\left(F^{n}\right){}^{k}(z)\underset{k\rightarrow\infty}{\rightarrow}z_{0}
\]
The set of ghost-periodic points of $F$ is $\PerG(F)$.
\end{defn}
\begin{rem}
By definition, $\PerG(F)\subset\omega(F)\cap(\mathcal{B}(F)-\alpha(F))$.
\end{rem}
We have a complete classification of the periodic regular components
of PMTs.
\begin{thm}
Let $U$ be a periodic regular component of period $n$ of the PMT
$F$. Then, only one of the following happens:
\begin{itemize}
\item Immediate basin of attraction, that is, exists an attracting periodic
point $z_{0}\in U$ such that for all $z\in U$ $\lim_{k\rightarrow\infty}(f^{n})^{k}(z)=z_{0}$.
\item Immediate parabolic basin, that is, exists a parabolic periodic point
$z_{0}\in\alpha(F)$ such that for all $z\in U$ $\lim_{k\rightarrow\infty}(f^{n})^{k}(z)=z_{0}$.
\item Immediate ghost-parabolic basin, that is, exists a ghost-periodic
point $z_{0}\in\partial U$ such that for all $z\in U$ $\lim_{k\rightarrow\infty}(f^{n})^{k}(z)=z_{0}$.
\item Rotation domain, that is, $F^{n}$ is an elliptic Möbius transformation
in $U$.
\item Neutral domain, that is, $F^{n}$ is the identity in $U$.
\end{itemize}
\end{thm}
\begin{rem}
In \cite{LerSie2019} the concepts of parabolic basin and ghost-parabolic
basin were not differentiated, but now we consider that it is important
to distinguish them due to their different dynamic behaviors.
\end{rem}
To conclude this Section, it is worth mentioning that there are examples
of PMT with wandering domains, with regular components of any connectivity,
with any number of regular components, or with pre-discontinuity set
of positive area (including the case of the whole sphere as pre-discontinuity
set), as discussed in \cite{LerSie2019}.

\section{Hyperbolicity and Expansivity\protect\label{sec:Hyperbolicity-and-Expansivity}}

It is well known that hyperbolic and structurally stable maps are
closely related, or most likely equivalent in the case of rational
maps. In this Section, we define and explore the concepts of hyperbolic
PMTs, to uncover their connections with structural stability.

Hyperbolic rational maps on $\Chat$ have only attracting and repelling
periodic points, and every periodic Fatou component is an immediate
attracting basin. The equivalent notion for PMTs can be defined using
this feature.
\begin{defn}
A PMT $F$ is \emph{hyperbolic} if $\PerA(F)\neq\emptyset$,\emph{
}$\PerN(F)=\emptyset$, \linebreak{}
$\PerG(F)=\emptyset$ and there are no wandering regular components.
\end{defn}
\begin{rem}
Note that the definition of hyperbolic PMT implies that every periodic
regular component is an immediate attracting basin.
\end{rem}
\begin{rem}
Prohibiting the existence of wandering components in the definition
of a hyperbolic PMT is essential, as these can lead to non-hyperbolic
dynamic behaviors. It is known that affine interval exchange transformations
(abbr. AIET) with wandering components (see \cite{BreEtAl2010,GutEtAl2008}),
where the component transformations are all contracting or expanding.
Let us construct a PMT $F$ as an extension of such AIET on $[0,1]$
to $\C$: take open discs $R_{k}$ with the corresponding interval
of the partition of the AIET as its diameter and an expanding transformation
$f$ on the exterior of the discs such that $f^{-1}(R_{k})\subset R_{1}$
for each $k$ and where $R_{1}$ is the element of the partition such
that $0\in\overline{R_{1}}.$ This PMT satisfies $\PerA(F)\neq\emptyset$
(at least $\infty\in\PerA(F)$),\emph{ }$\PerN(F)=\emptyset,$ and
$\PerG(F)=\emptyset,$ but the wandering components accumulate in
$\alpha(F)$. Therefore, there exist $z\in\Rc(F)$ such that their
orbits do not converge to a periodic attracting point.
\end{rem}
Unlike hyperbolic rational maps, hyperbolic PMTs might not have repelling
periodic points.
\begin{example}
Let
\[
F(z)=\begin{cases}
\lambda z & \mathrm{if}\,z\in\mathbb{D}\\
\frac{1}{\lambda}z & \mathrm{if}\,z\in\Chat-\overline{\mathbb{D}},
\end{cases}
\]
where $\mathbb{D}=\left\{ z\in\mathbb{C}:\,|z|<1\right\} $ and $\lambda\in\mathbb{D}-\left\{ 0\right\} $.

Then $0$ and $\infty$ are attracting fixed points with $\mathbb{D}$
and $\Chat-\overline{\mathbb{D}}$ as attracting basins, respectively.
Since $\Bc(F)=B(F)=\partial\mathbb{D}$, there are no repelling or
neutral periodic points. That is, $F$ is a hyperbolic PMT without
repelling periodic points.

\rule[1ex]{0.1\paperwidth}{0.25pt}
\end{example}
The hyperbolic behavior in PMTs arises from the loxodromic component
functions. However, not all component functions need to be loxodromic
for the PMTs to exhibit hyperbolic characteristics, as demonstrated
in the following example.
\begin{example}
Let
\[
F(z)=\begin{cases}
z+2 & \mathrm{if}\,z\in\mathbb{D}\\
2z & \mathrm{if}\,z\in\Chat-\overline{\mathbb{D}}.
\end{cases}
\]

Then $\Rc(F)=\mathbb{D}\cup(\Chat-\overline{\mathbb{D}})$. The only
periodic component is $\Chat-\overline{\mathbb{D}}$, the immediate
attracting basin of $\infty$ the unique attracting fixed point of
$F$. The regular component $\mathbb{D}$ is preperiodic. The transformation
$z\mapsto z+2$ is not loxodromic, but $F$ is clearly hyperbolic.

\rule[1ex]{0.1\paperwidth}{0.25pt}
\end{example}
On the other hand, a PMT with all its component functions loxodromic
is not necessarily hyperbolic.
\begin{example}
Let
\[
F(z)=\begin{cases}
\frac{1}{2}z & \mathrm{if}\,z\in R_{1}\\
2z & \mathrm{if}\,z\in R_{2},
\end{cases}
\]
where $R_{1}=\left\{ z\in\mathbb{C}:\,|z-1|<1\right\} $ and $R_{2}=\Chat-\overline{R_{1}}$.

We have that $\Rc(F)=R_{1}\cup R_{2}$. Note that both component functions
are loxodromic, but $0$ is a ghost-periodic point. Therefore, $F$
is not hyperbolic.

\rule[1ex]{0.1\paperwidth}{0.25pt}
\end{example}
For hyperbolic rational maps on $\Chat$, the dynamical behavior can
be associated with certain conditions regarding the post-critical
set. PMTs do not have critical points; however, the dynamical behavior
can be associated with the $\omega$-limit set.
\begin{thm}
\label{thrm:HyperOmega}Let $F$ be a PMT. Then the following conditions
are equivalent:
\begin{enumerate}
\item $F$ is hyperbolic.
\item $\omega(F)=\PerA(F)\neq\emptyset$.
\end{enumerate}
\end{thm}
\begin{proof}
~
\begin{enumerate}
\item Let $F$ be hyperbolic. By the definition of $\omega-$limit we have
$\omega(F)=\PerA(F)\cup\PerE(F)\cup\PerID(F)\cup\PerP(F)\cup\PerG(F)$.
Thus, we conclude $\omega(F)=\PerA(F)\neq\emptyset$.
\item Suppose that $\omega(F)=\PerA(F)\neq\emptyset$. By the definitions
of ghost-periodic point and $\omega$-limit set, $\PerP(F)=\PerE(F)=\PerID(F)=\PerG(F)=\emptyset$.
That is, $F$ is hyperbolic.
\end{enumerate}
\end{proof}
\begin{rem}
Note that if $F$ is a hyperbolic PMT, by the incise (2) of Theorem
\ref{thrm:HyperOmega}, we have $\omega(F)\cap\Bc(F)=\emptyset$ because
there are no parabolic periodic points, no ghost-periodic points,
and no wandering components.
\end{rem}
Contrary to the conjectured equivalence between being hyperbolic and
structurally stable in rational maps on the Riemann sphere, there
exist hyperbolic PMTs which are not structurally stable.
\begin{example}
\label{example:HiperNoSS}Let
\[
F_{\lambda}(z)=\begin{cases}
f_{1}(z) & \mathrm{if}\,z\in R_{1}\\
f_{2}(z) & \mathrm{if}\,z\in R_{2},
\end{cases}
\]
where $f_{1}(z)=\lambda z+\lambda$, $f_{2}(z)=\frac{6i\lambda z-1}{z+6i\lambda}$,
$R_{1}=\left\{ z:\,|1-z|<1\right\} $ and $R_{2}=\Chat-\overline{R_{1}}$.
$f_{1}$ and $f_{2}$ are both loxodromic when $0<|\lambda|<1$.

Let $\lambda_{0}=\frac{1}{2}$. Then, there exists a neighborhood
$\Nc_{\lambda_{0}}\subset\C$ such that $f_{1}$ and $f_{2}$ are
loxodromic. The fixed points of $f_{1}$ are $z_{\lambda}=\frac{\lambda}{1-\lambda}$
(attracting) and $\infty$ (repelling), while the fixed points of
$f_{2}$ are always $i$ (attracting) and $-i$ (repelling). Then
the neighborhood $\Nc_{\lambda_{0}}$ can be adjusted in such a way
that $z_{\lambda}\in R_{1}$ for all $\lambda\in\Nc_{\lambda_{0}}$.
Therefore, $R_{1}$ must contain an immediate basin of attraction
for the fixed point $z_{\lambda}$. Even more, for all $\lambda\in\Nc_{\lambda_{0}}$
we have $i,-i\in R_{2}$, implying that $R_{2}$ contains an immediate
basin of attraction for the fixed point $i$, and that $-i\in\alpha(F)$.

\begin{figure}[h]
\begin{centering}
\includegraphics{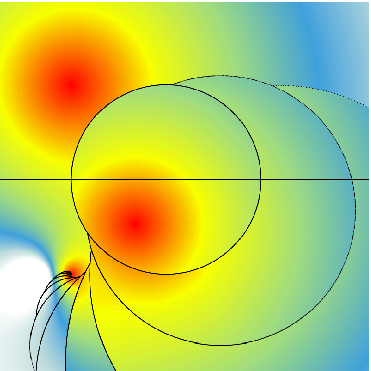}\hspace*{0.25cm}\includegraphics{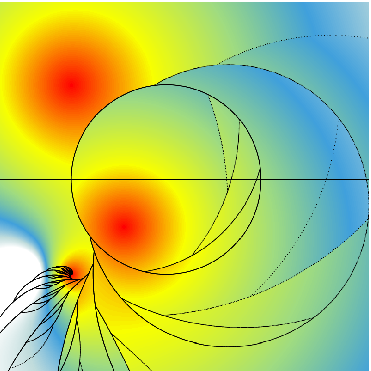}
\par\end{centering}
\caption{\protect\label{fig:Fig1}The pre-discontinuity and regular sets of
$F_{\lambda}$ described in Example \ref{example:HiperNoSS}.\protect \\
Left: With $\lambda=\frac{1}{2}-0.223i$. $R_{1}$ is the immediate
basin of attraction of $z_{\lambda}$. Right: With $\lambda=\frac{1}{2}-(0.223+\varepsilon)i$,
$0<\varepsilon\ll1$. $R_{1}$ contains several regular components.}
\end{figure}

For each $\lambda\in\Nc_{\lambda_{0}}$, let $A_{\lambda}$ be the
immediate basin of attraction of $z_{\lambda}\in R_{1}$, $U_{\lambda}=\bigcup_{n\geq0}F^{-n}(A_{\lambda})$
and $V_{\lambda}=\Rc(F)-U_{\lambda}$. Then, $F^{n}(z)\underset{n\rightarrow\infty}{\rightarrow}z_{\lambda}$
for all $z\in U_{\lambda}$ and $F^{n}(z)\underset{n\rightarrow\infty}{\rightarrow}i$
for all $z\in V_{\lambda}$. Therefore, $F_{\lambda}$ has only three
periodic points, all of which are fixed: $z_{\lambda}$, $i$ and
$-i$. Furthermore, these fixed points are either attracting or repelling,
so $F_{\lambda}$ is hyperbolic.

On the other hand, varying $\lambda$ within $\Nc_{\lambda_{0}}$,
there exist maps such that the immediate basin of attraction of $z_{\lambda}$
is exactly $R_{1}$, as well as maps where $R_{1}$ contains several
regular components. Clearly, these maps cannot be conjugated. Then,
there exists parameters $\lambda'\in\Nc_{\lambda_{0}}$ where the
the aforementioned bifurcation occurs, indicating that $F_{\lambda}$
is not structurally stable in neighborhoods $\Nc_{\lambda'}\subset\Nc_{\lambda_{0}}$.

To clarify this example, Figure \ref{fig:Fig1} illustrates the approximations
of the pre-discontinuity sets of $F_{\lambda}$ in black, and the
attracting fixed points $z_{\lambda}\in R_{1}$ and $i$, as well
the repelling fixed point $-i\in\alpha(F)$, in the center of the
red spots.

\rule[1ex]{0.1\paperwidth}{0.25pt}
\end{example}
For PMTs, a similar definition to expanding rational maps exists,
utilizing points in the pre-discontinuity set where all iterations
of the map are defined and differentiable.
\begin{defn}
A PMT $F$ is \emph{$\alpha$-expanding} if there is a $N\geq1$ such
that $|(F^{N})'(z)|_{s}>1$ (where $|\cdot|_{s}$ is the normalized
spherical norm) for all $z\in\alpha(F)$.
\end{defn}
In contrast to rational maps on the Riemann sphere, the properties
of being hyperbolic and $\alpha$-expanding are not equivalent for
PMTs, as demonstrated in the following examples.
\begin{example}
Let
\[
F(z)=\begin{cases}
\lambda z & \mathrm{if}\,z\in\mathbb{D}\\
\frac{1}{\lambda}z & \mathrm{if}\,z\in\Chat-\overline{\mathbb{D}},
\end{cases}
\]
where $\lambda\in\mathbb{D}-\left\{ 0\right\} $.

As observed earlier, $F$ is hyperbolic and is not $\alpha$-expanding
since $\alpha(F)=\emptyset$.

\rule[1ex]{0.1\paperwidth}{0.25pt}
\end{example}
\begin{example}
\label{example:ExpandNoHyper}There exists $\alpha$-expanding but
non-hyperbolic PMT; this is due to the absence of incompatibility
between being expanding and having elliptic, identity, and ghost-periodic
points.

\begin{figure}[h]
\begin{centering}
\includegraphics{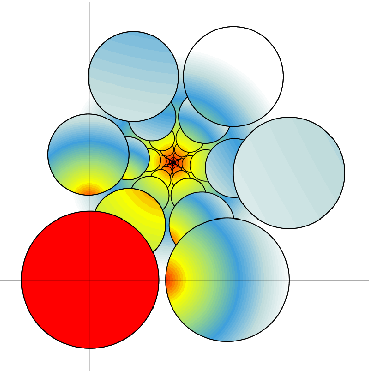}
\par\end{centering}
\caption{\protect\label{fig:Fig2}Pre-discontinuity set (depicted in black)
and regular set (drawn with colors) of $F$ from example \ref{example:ExpandNoHyper}.}
\end{figure}
\end{example}
Let
\[
F(z)=\begin{cases}
e^{\frac{2}{3}\pi i}z & \mathrm{if}\,z\in R_{1}\\
\frac{10}{9}e^{\frac{2}{3}\pi i}(1-z) & \mathrm{if}\,z\in R_{2},
\end{cases}
\]

where $R_{1}=\left\{ z:\,|z|<\frac{1}{2}\right\} $ and $R_{2}=\Chat-\overline{R_{1}}$.

$R_{1}$ is a rotation domain where $0$ is an elliptic fixed point,
and $z_{0}=\frac{\lambda}{\lambda+1}$, with $\lambda=\frac{10}{9}e^{\frac{2}{3}\pi i}$,
is a repelling fixed point.

Clearly, $\alpha(F)=\left\{ z_{0}\right\} $ and thus $F$ is \emph{$\alpha$-}expanding
but not hyperbolic.

\rule[1ex]{0.1\paperwidth}{0.25pt}

In the case of non-hyperbolic and non $\alpha$-expansive PMTs, strange
dynamic behaviors may occur, as illustrated in the following example.
\begin{example}
There exist a non $\alpha$-expanding PMTs with two repelling fixed
points and forward invariant subsets $A\subset\alpha(F)$ such that
$F|_{A}$ is conjugated to an irrational rotation. This map has no
regular components.

For the PMT
\[
F(z)=\begin{cases}
2z & \mathrm{if}\,z\in\left\{ z:\,|z|<1\right\} \\
\frac{2}{3}z & \mathrm{if}\,z\in\left\{ z:\,|z|>1\right\} 
\end{cases}
\]
it has been proven that $F|_{[\frac{2}{3},2)}$ is topologically conjugated
with an irrational rotation in $S^{1}$ and $F$ behaves the same
in all rays from $0$ to $\infty$ (see \cite{LerSie2019}).

Therefore, for all $z\in\left\{ z\in\C\,:\,\frac{2}{3}\leq|z|\leq2\right\} \cap\alpha(F)$,
there cannot exist $N\geq0$ such that $|F^{N}(z)|_{s}>1$ since $F|_{O(z,F)}$
is conjugated with an irrational rotation on an orbit subset of $S^{1}$.

On the other hand, $\Per(F)=\Fix(F)=\left\{ 0,\infty\right\} $ are
repelling.

\rule[1ex]{0.1\paperwidth}{0.25pt}
\end{example}
As has been exposed, there is a non-equivalence between hyperbolic
and \emph{$\alpha$-}expanding notions for PMTs; therefore, they cannot
be studied as a single concept. The possibility for generating drastic
changes in the regular set through perturbations of hyperbolic maps
renders an equivalence of this notion with structural stability impossible.
Lastly, the compatibility between the existence of elliptic, identity,
and ghost-periodic points and the property of being \emph{$\alpha$-}expanding,
implies that such maps are not necessarily structurally stable.

\section{Parameter space of PMTs and conjugations}

The parameter space of PMTs $F=(\left\{ R_{k}\right\} _{k=1}^{K},\,\left\{ f_{k}\right\} _{k=1}^{K})$
is determined by the maps $F\vert_{R_{k}}=f_{k}\in PSL(2,\C)$ and
the elements $R_{k}$ of the partition in $\Chat$. For the partition,
it suffices to consider the space of discontinuity sets $B=\bigcup_{k=1}^{K}\partial R_{k}$
as compact subsets of $\Chat$. Thus, we can establish the following
\begin{defn}
The parameter space of PMTs over a partition of $\Chat$ in $K>1$
parts is
\[
X_{PMT,K}=\overbrace{PSL(2,\C)\times\dots\times PSL(2,\C)}^{K\,\mathrm{times}}\times\Pc_{K}(\Chat)
\]
with the product topology, where $\Pc_{K}(\Chat)$ is the space of
the discontinuity sets whose associated partitions in $\Chat$ has
$K$ parts.
\end{defn}
\begin{rem}
$\Pc_{K}(\Chat)$ is a subset of the space of non-empty compact subsets
of $\Chat$, equipped with the Hausdorff metric. However, $\Pc_{K}(\Chat)$
can also be regarded as a Teichmüller space since each $B\in\Pc_{K}(\Chat)$
defines a set of regions $R_{k}$ which are Riemann surfaces, consequently
$\Pc_{K}(\Chat)\subset Teich(R_{1})\times Teich(R_{2})\times\dots\times Teich(R_{K})$.
Furthermore, $\Pc_{K}(\Chat)$ is a complex manifold because every
$R_{k}$ is a hyperbolic Riemann surface, as follows from the Bers
embedding theorem (see for example \cite{GarLak2000}). In this work,
the holomorphic structure of this parameter space will be particularly
beneficial to us.
\end{rem}
As usual, $F,G\in X_{PMT,K}$ are topologically conjugated if there
exists a homeomorphism $h:\Chat\,\selfarrow$ such that $h\circ F=G\circ h$.
The next result follows immediately.
\begin{thm}
If $F,G\in X_{PMT,K}$ are topologically conjugated by a homeomorphism
$h:\Chat\,\selfarrow$, then $B(G)=h(B(F))$, $\Bc(G)=h(\Bc(F))$,
$\alpha(G)=h(\alpha(F))$, and $\Rc(G)=h(\Rc(F))$.
\end{thm}

\section{Structural Stability in $PSL(2,\mathbb{C})^{K}$}

In this Section, we will examine the stability of all PMTs by fixing
the discontinuity set $B$ and perturbing the component functions.
The corresponding parameter space with this fixture is $PSL(2,\C)^{K}\cong PSL(2,\C)^{K}\times\left\{ B\right\} \subset X_{PCM,K}$.

Now, we can establish the following
\begin{defn}
A PMT $F=\left(\left\{ R_{k}\right\} _{k=1}^{K},\left\{ f_{k}\right\} _{k=1}^{K}\right)$
is \emph{structurally stable in $PSL(2,\C)^{K}$} if there exists
a neighborhood $\mathcal{N}_{(f_{1},\dots,f_{K})}\subset PSL(2,\C)^{K}$
such that for every element $(g_{1},\dots,g_{K})\in\mathcal{N}_{(f_{1},\dots,f_{K})}$,
there exists a homeomorphism $h:\Chat\rightarrow\Chat$ such that
$h\circ F=G\circ h$ in the conformality region $R(F)$, and the discontinuity
set is fixed (that is $B(F)=B(G)$), where $G$ is the corresponding
PMT $\left\{ \left\{ R_{k}\right\} _{k=1}^{K},\left\{ g_{k}\right\} _{k=1}^{K}\right\} $.
\end{defn}
One of the results in \cite{LerSie2019} establishes the sufficiency
of the structural stability in $PSL(2,\C)^{K}$ if $\left\langle f_{1},\dots,f_{K}\right\rangle $
is a structurally stable group and the boundary set is contained within
a fundamental region of that group. However, the structural stability
of PMTs can indeed be achieved without any additional requirements
on the group $\left\langle f_{1},\dots,f_{K}\right\rangle $, utilizing
several strong hypotheses as outlined below. 
\begin{thm}
\label{thrm:StrucStabPSL2CK}Let $F=\left(\left\{ R_{k}\right\} _{k=1}^{K},\left\{ f_{k}\right\} _{k=1}^{K}\right)$
be a PMT such that
\begin{enumerate}
\item each component transformation $f_{k}$ is loxodromic,
\item $F$ is hyperbolic,
\item for each $k$, one of the following statements holds
\begin{enumerate}
\item $f_{k}^{-1}(B(F))\cap R_{k}=f_{k}^{-1}(B(F))$,
\item $f_{k}^{-1}(B(F))\cap R_{k}=f_{k}^{-1}(B_{j})$ for some connected
component $B_{j}$ of $B(F)$, or
\item $f_{k}^{-1}(B(F))\cap R_{k}=\emptyset$.
\end{enumerate}
\item for all $n>1$ and for each connected component $C_{i}$ of $F^{-n}(B(F))$,
$F^{n}(C_{i})=B_{j}$ for some connected component $B_{j}$ of $B(F)$,
being $F^{n}|_{C_{i}}$ a Möbius transformation;
\end{enumerate}
then $F$ is structurally stable in $PSL(2,\C)^{K}$.
\end{thm}
\begin{rem}
Hypothesis (1) is essential since parabolic and elliptic Möbius transformations
are not structurally stable. Hypothesis (2) is evidently necessary,
as discussed in Section \ref{sec:Hyperbolicity-and-Expansivity}.
Hypotheses (3) and (4) establish a Schottky-like behavior for the
PMT. The hypothesis (3) is precisely the hypothesis (4) for the case
$n=1$, but they are presented separately for clarity.
\end{rem}
\begin{proof}
Small perturbations of loxodromic maps remain loxodromic, and for
this reason, hyperbolic PMTs with loxodromic component functions continue
to be hyperbolic. The action of $PSL(2,\C)$ in $\Chat$ is continuous,
ensuring that disjoint subsets remain disjoint under the action of
maps in a small neighborhood of the component functions. Thus, hypotheses
(1), (2), (3), and (4) permit us to take a neighborhood $\Nc_{F}=\Nc_{(f_{1},\dots,f_{K})}\subset PSL(2,\C)^{K}$
such that for all $(g_{1},\dots,g_{K})\in\Nc_{F}$, the defined PMT
$G\equiv\left\{ \left\{ g_{k}\right\} _{k=1}^{K},\left\{ R_{k}\right\} _{k=1}^{K}\right\} $
also satisfy hypotheses (1), (2), (3), and (4). 

Let $B=B(F)$. We construct $\varphi:\mathcal{N}_{F}\times E\rightarrow\Chat$,
a holomorphic motion of $E=\left(\bigcup_{n\geq0}F^{-n}(B)\right)\cup\PerA(F)$
as follows. For $\lambda=(g_{1},\dots,g_{K})\in\mathcal{N}_{F}$ with
associated PMT $G$ and $z\in E$, define 
\[
\varphi(\lambda,z)=\begin{cases}
z & \mathrm{if}\,z\in B\\
G^{-n}\circ F^{n}(z) & \mathrm{if}\,z\in F^{-n}(B),\,n>0\\
w_{z} & \mathrm{if}\,z\in\PerA(F)
\end{cases}
\]
where $G^{-n}\circ F^{n}$ is a composition $g_{k_{1}}^{-1}\circ\dots\circ g_{k_{n}}^{-1}\circ f_{k_{n}}\circ\dots\circ f_{k_{1}}$
and $w_{z}$ is the attracting fixed point of $G^{n}=g_{k_{n}}\circ\dots\circ g_{k_{1}}$
associated to the corresponding attracting fixed point $z$ of $F^{n}=f_{k_{n}}\circ\dots\circ f_{k_{1}}$.

Observe that if $z\in F^{-n}(B)$, then $F^{n}(z)\in B$ and $G^{-n}\circ F^{n}(z)\in G^{-n}(B)\subset\Bc(G)$.
Using hypotheses (3) and (4) each function $\varphi_{\lambda}=\varphi(\lambda,\_)$
is an injection on $\Chat$ because $\varphi_{\lambda}$ is defined
by a single Möbius transformation in each set homeomorphic to $B$
or to $B_{j}$ (component of $B$) forming $F^{-n}(B)$, or is the
identity in $B$, or is the bijection between attracting periodic
points. Such a bijection between attracting periodic points is possible
due to the hypotheses, since $F$ and $G$ do not have parabolic,
elliptic or identity periodic points, and regular components are preserved.

The function $\lambda\mapsto\varphi(\lambda,z)$ is a composition
of the Möbius transformations $g_{1}^{-1},\dots$, $g_{K}^{-1}$,
$f_{1},\dots$, $f_{K}$ with the parameters moving holomorphically;
thus, $\varphi(\lambda,z_{0})$ is a holomorphic function of $\lambda$
for each $z_{0}\in E$. If $\lambda_{0}$ is the element associated
to $F$, is clear that $\varphi(\lambda_{0},z)=z$.

Using the Bers-Royden extension theorem (see \cite{BerRoy1986}),
$\varphi$ has an extension to a holomorphic motion $\Phi$ of $\Chat$.
This can be accomplished as follows:
\begin{itemize}
\item First, restrict $\varphi$ to a disc $D\subset\mathcal{N}_{F}$, then
transforming $D$ to $\mathbb{D}$ with via an affine map, and finally
restrict to $D(0,\frac{1}{3})=\left\{ z:\,|z|<\frac{1}{3}\right\} $.
Thus, $\varphi|_{D(0,\frac{1}{3})\times E}$ can be extended to $\Phi:D(0,\frac{1}{3})\times\Chat\rightarrow\Chat$,
as stated in the Bers-Royden extension theorem.
\item Furthermore, for each $\lambda\in D(0,\frac{1}{3})$, the map $z\mapsto\Phi(\lambda,z)$
is a quasi-conformal homeomorphism $h_{\lambda}:\Chat\rightarrow\Chat$,
which can be uniquely chosen such that the Beltrami differential $\mu(h_{\lambda})$
is harmonic in $\Chat-\overline{E}$. By the connectivity of $\mathcal{N}_{F}$
and the uniqueness of $\Phi$, the holomorphic motion $\Phi$ can
be adapted and extended to $\mathcal{N}_{F}\times\Chat$.
\end{itemize}
By construction, $h_{\lambda}=\Phi(\lambda,\_)$ conjugates $F$ with
$G$:
\begin{itemize}
\item If $z\in B$, then $F$ and $G$ are undefined on $z$. Since $h_{\lambda}|_{B}\equiv Id|_{B}$,
it follows that $h_{\lambda}\circ F$ and $G\circ h_{\lambda}$ are
undefined on $z$.
\item If $z\in F^{-1}(B)$, then $F(z)\in B$. By definition of $h_{\lambda}$
using $\varphi$:
\begin{itemize}
\item $h_{\lambda}\circ F(z)=F(z)$.
\item $G\circ h_{\lambda}(z)=G\circ G^{-1}\circ F(z)=F(z)$.
\end{itemize}
\item If $z\in F^{-n}(B)$ for some $n>1$, then $F(z)\in F^{-n+1}(B)$.
By the definition of $h_{\lambda}$ using $\varphi$:
\begin{itemize}
\item $h_{\lambda}\circ F(z)=G^{-n+1}\circ F^{n-1}(F(z))=G^{-n+1}\circ F^{n}(z)$. 
\item $G\circ h_{\lambda}(z)=G\circ G^{-n}\circ F^{n}(z)=G^{-n+1}\circ F^{n}(z)$.
\end{itemize}
\item If $z\in\PerA(F)$, is it the attracting fixed point of some composition
$f_{k_{n}}\circ\dots\circ f_{k_{1}}$. Observe that $z=F^{n}(z)\in R_{k_{1}}$,
then
\begin{enumerate}
\item $F(z)$ is the attracting fixed point of $f_{k_{1}}\circ f_{k_{n}}\circ\dots\circ f_{k_{2}}$.
Therefore, $h_{\lambda}\circ F(z)$ is the attracting fixed point
of $g_{k_{1}}\circ g_{k_{n}}\circ\dots\circ g_{k_{2}}$.
\item $h_{\lambda}(z)$ is the attracting fixed point of $g_{k_{n}}\circ\dots\circ g_{k_{1}}$
and $h_{\lambda}(z)=G^{n}(h_{\lambda}(z))\in R_{k_{1}}$. Therefore,
$G\circ h_{\lambda}(z)$ is the attracting fixed point of $g_{k_{1}}\circ g_{k_{n}}\circ\dots\circ g_{k_{2}}$.
\end{enumerate}
\item The function in $\mathcal{N}_{F}\times\Chat$ given by
\[
\tilde{h}_{\lambda}(z)=\begin{cases}
g_{k}^{-1}\circ h_{\lambda}\circ f_{k} & \mathrm{if}\,z\in R_{k}\\
z & \mathrm{if}\,z\in B
\end{cases}
\]
is also an extension of the holomorphic motion $\varphi$, with harmonic
Beltrami differential since $f_{k}$ and $g_{k}^{-1}$ are holomorphic.
By the uniqueness of the Bers-Royden extension under such condition,
we have $\tilde{h}_{\lambda}=h_{\lambda}$.

Therefore, if $z\in\Chat-E$, then $z\in R_{k}$ for some $k$, and
we can conclude
\[
h_{\lambda}\circ F(z)=h_{\lambda}\circ f_{k}(z)=g_{k}\circ h_{\lambda}(z)=G\circ h_{\lambda}(z).
\]

\end{itemize}
\end{proof}
\begin{rem}
This theorem and its proof serve as the foundation and inspiration
of the statement and proof of the final Theorem \ref{thrm:GenStrucStab},
concerning structural stability in the general case of the parameter
space $X_{PMT,K}$.
\end{rem}
\begin{example}
\label{example:SS}Let
\[
F(z)=\begin{cases}
f_{1}(z) & \mathrm{if}\,z\in R_{1}=\left\{ z:\,|z|<\frac{2}{5}\right\} \\
f_{2}(z) & \mathrm{if}\,z\in R_{2}=\Chat-\overline{R_{2}},
\end{cases}
\]
where $f_{1}(z)=\frac{(1+i)z+i}{-iz+(1-i)}$ and $f_{2}(z)=\frac{(1+i)z-i}{iz+(1-i)}$
are loxodromic maps. $1$ is a parabolic fixed point for $F$, thus
$F$ does not satisfy the hyperbolicity hypothesis of the previous
theorem.

On the other hand, $f_{1}$ and $f_{2}$ can be slightly perturbed
to remain loxodromic maps, ensuring that the corresponding PMT contains
only attracting and repelling fixed points, with no parabolic fixed
points. These perturbations can be executed in a manner that fulfills
the hypotheses of the previous theorem, making all these perturbed
PMTs structurally stable in $PSL(2,\C)^{2}$. Such perturbed PMTs
exhibit the following dynamic characteristics:
\begin{itemize}
\item $\Bc(F)$ consists of the union of an infinite number of disjoint
circles, along with the $\alpha$-limit set.
\item They possess a single attracting fixed point and a single repelling
fixed point, both located at the centers of the spots colored in red.
\item They feature a unique immediate basin of attraction: the exterior
of the discs whose boundaries form $\Bc(F)$.
\item The regular components, which are the interiors of the discs forming
$\Bc(F)$, are pre-periodic.
\end{itemize}
\end{example}
\begin{figure}[h]
\begin{centering}
\includegraphics[scale=0.9]{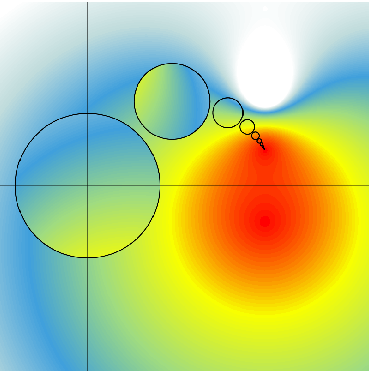}\hspace*{0.25cm}\includegraphics[scale=0.9]{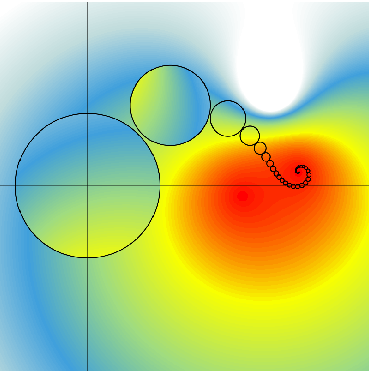}
\par\end{centering}
\begin{centering}
\includegraphics[scale=0.9]{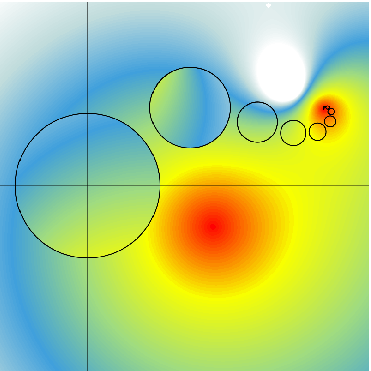}\hspace*{0.25cm}\includegraphics[scale=0.9]{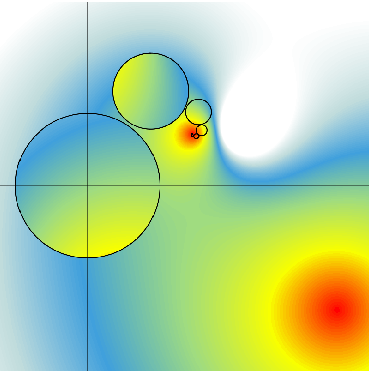}
\par\end{centering}
\caption{\protect\label{fig:Fig3}The pre-discontinuity and regular sets of
$F$ composed of $f_{1}$ and $f_{2}$, from Example \ref{example:SS}.\protect \\
Top left: $f_{1}(z)=\frac{(1+i)z+1.02i}{-1.02iz+(1-i)}$ and $f_{2}(z)=\frac{(1+i)z-1.02i}{1.02iz+(1-i)}$.\protect \\
Top right: $f_{1}(z)=\frac{(1+i)z+c}{-cz+(1-i)}$ and $f_{2}(z)=\frac{(1+i)z-c}{cz+(1-i)}$,
with $c=0.99+0.01i$.\protect \\
Bottom left: $f_{1}(z)=\frac{(1+i)z+i}{-iz+(1-i)}$ and $f_{2}(z)=\frac{(0.8+i)z-i}{iz+(1-i)}$.\protect \\
Bottom right: $f_{1}(z)=\frac{(1+i)z+i}{-iz+(1-i)}$ and $f_{2}(z)=\frac{(1.1+i)z+0.1-i}{(-0.1+i)z+(0.9-i)}$.}
\end{figure}
In the images from Figure \ref{fig:Fig3}, the approximations of the
pre-discontinuity sets of perturbations of $F$ are depicted in black.

\rule[1ex]{0.1\paperwidth}{0.25pt}

\section{$\Bc$-Stability}

Before the study of general structural stability of PMTs, we will
define and analyze a type of stability analogous to the J-stability
of rational maps.

First, we define holomorphic families of PMTs, where the corresponding
parameter space is necessarily a complex manifold.
\begin{defn}
A family of PMTs $\left\{ F_{\mu,\lambda}:\Chat\,\selfarrow\right\} _{(\mu,\lambda)\in Y\times X}$,
parameterized by $(\mu,\lambda)\in Y\times X$ where $Y$ and $X$
are complex manifolds, is a \emph{holomorphic family} if
\begin{itemize}
\item There exists a holomorphic motion of the discontinuity set $B(F_{\mu_{0},\lambda})\in\Pc_{K}(\Chat)$,
parameterized by $(Y,\mu_{0})$ over the discontinuity sets of $F_{\mu,\lambda}$.
\item The function $Y\times X\times R(F_{\mu,\lambda})\rightarrow R(F_{\mu,\lambda})$,
given by $(\mu,\lambda,z)\mapsto F_{\mu,\lambda}(z)$, is holomorphic.
\end{itemize}
\end{defn}
In a manner analogous to the definition of holomorphic motion for
Julia sets, it can also be defined for the pre-discontinuity sets
of PMTs.
\begin{defn}
Given a holomorphic family of PMTs $\left\{ F_{\mu,\lambda}:\Chat\,\selfarrow\right\} _{(\mu,\lambda)\in Y\times X}$,
the pre-discontinuity sets $\Bc(F_{\mu,\lambda})$ \emph{move holomorphically
}if there exists a holomorphic motion
\[
\left\{ \varphi_{\mu,\lambda}:\Bc(F_{\mu_{0},\lambda_{0}})\rightarrow\Chat\right\} _{(\mu,\lambda)\in Y\times X}
\]
 such that 
\[
\varphi_{\mu,\lambda}\left(\Bc(F_{\mu_{0},\lambda_{0}})\right)=\Bc(F_{\mu,\lambda}),
\]
\[
\varphi_{\mu,\lambda}\circ F_{\mu_{0},\lambda_{0}}|_{\Bc(F_{\mu_{0},\lambda_{0}})-B(F_{\mu_{0},\lambda_{0}})}=F_{\mu,\lambda}\circ\varphi_{\mu,\lambda}|_{\Bc(F_{\mu_{0},\lambda_{0}})-B(F_{\mu_{0},\lambda_{0}})},
\]
and 
\[
\varphi_{\mu,\lambda}(B(F_{\mu_{0},\lambda_{0}}))=B(F_{\mu,\lambda}).
\]
The pre-discontinuity sets $\Bc(F_{\mu,\lambda})$ \emph{move holomorphically
at} $(\mu_{0},\lambda_{0})$ if they move holomorphically in some
neighborhood $\Nc_{(\mu_{0},\lambda_{0})}\subset Y\times X$.
\end{defn}
\begin{rem}
Note that the holomorphic motion $\varphi_{\mu,\lambda}$ may not
respect the dynamics across the entire set $\Bc(F_{\mu,\lambda})$,
due to the undefinition of $F_{\mu,\lambda}$ in $B(F_{\mu,\lambda})$.
\end{rem}
Now, we can be defined the concept of $\Bc$-stability.
\begin{defn}
A PMT $F$ is \emph{$\Bc$-stable} if $\Bc(F)$ moves holomorphically.
\end{defn}
As anticipated, there exist PMTs that are $\Bc$-stable but not structurally
stable, as demonstrated below.
\begin{example}
\label{example:HoloMot}Let
\[
F_{\mu,\lambda}(z)=\begin{cases}
f_{1}(z) & \mathrm{if}\,z\in R_{1}\\
f_{2}(z) & \mathrm{if}\,z\in R_{2},
\end{cases}
\]
where $f_{1}(z)=\frac{(1+i)z+\lambda}{-\lambda z+(1-i)}$, $f_{2}(z)=\frac{(1+i)z-\lambda}{\lambda z+(1-i)}$,
$R_{1}=\left\{ z:\,|z-\mu|<\frac{1}{3}\right\} $ and $R_{2}=\Chat-\overline{R_{1}}$,
with $(\mu,\lambda)\in\left\{ \lambda:\,|\lambda|<\frac{1}{10}\right\} \times\left\{ \mu:\,|\mu-i|<\frac{1}{10}\right\} =Y\times X.$
Clearly, $F_{\mu,\lambda}$ is a holomorphic family of PMTs.
\begin{figure}[h]
\begin{centering}
\includegraphics{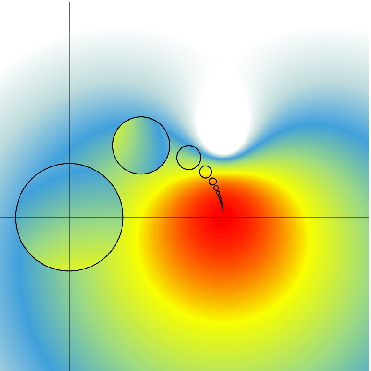}\hspace*{0.25cm}\includegraphics{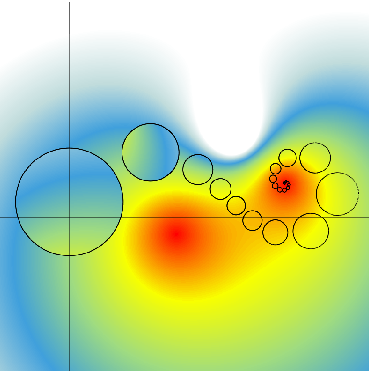}
\par\end{centering}
\caption{\protect\label{fig:Fig4}Holomorphic motion of $\Bc(F_{0,i})$, from
Example \ref{example:HoloMot}.\protect \\
Left: With $\mu=0$ and $\lambda=i$, $F_{0,i}$ has a unique fixed
point at $z=1$, which is parabolic. Right: With $\mu\approx0$ and
$\lambda\approx i$, $F_{\mu,\lambda}$ has two fixed points: $\frac{i+\sqrt{-1-\lambda{{}^2}}}{\lambda}$
attracting, and $\frac{i-\sqrt{-1-\lambda{{}^2}}}{\lambda}$ repelling.}
\end{figure}

A holomorphic motion $\varphi_{\mu,\lambda}:\Bc(F_{0,i})\rightarrow\Chat$
can be expressed as
\[
\varphi_{\mu,\lambda}(z)=\begin{cases}
z+\mu & z\in B(F_{0,i})\\
F_{\mu,\lambda}^{-N}(F_{0,i}^{N}(z)+\mu) & z\in\Bc(F_{0,i})-\alpha(F_{0,i})\\
\frac{i-\sqrt{-1-\lambda{{}^2}}}{\lambda} & z\in\alpha(F_{0,i})=\left\{ 1\right\} 
\end{cases}
\]
Then $\Bc(F_{0,i})$ moves holomorphically, but $F_{0,i}$ and $F_{\mu,\lambda}$
are not conjugated, for $(\mu,\lambda)$ as close to $(0,i)$ as we
desire.

In Figure \ref{fig:Fig4}, approximations of the pre-discontinuity
sets of $F_{\mu,\lambda}$ are depicted in black, with fixed points
located at the center of the red spots.

\rule[1ex]{0.1\paperwidth}{0.25pt}
\end{example}
\begin{rem}
From the previous example, we can observe that in the holomorphic
motions of PMTs, parabolic points can be transformed into repelling
points, unlike the holomorphic motions of rational maps.
\end{rem}
A consequence of the earlier definitions and the invariance of the
$\alpha$-limit set, is the following corollary.
\begin{cor}
If a PMT $F$ is $\Bc$-stable, then there exists a holomorphic motion
\[
\left\{ \varphi_{\mu,\lambda}:\alpha(F)\rightarrow\Chat\right\} _{(\mu,\lambda)\in\Nc\subset Y\times X}
\]
 such that $\varphi_{\mu,\lambda}\left(\alpha(F)\right)=\alpha(F_{\mu,\lambda})$
and 
\[
\varphi_{\mu,\lambda}\circ F|_{\alpha(F)}=F_{\mu,\lambda}\circ\varphi_{\mu,\lambda}|_{\alpha(F)}.
\]
\end{cor}
\begin{rem}
This corollary can be interpreted as follows: $\Bc$-stability implies
structural stability in the $\alpha$-limit set, as the corresponding
holomorphic motion preserves the dynamics of the $\alpha$-limit set.
\end{rem}
As usual, the concept of $\Bc$-stability across the entire parameter
space of PMTs is referred as $\Bc$-structural stability.
\begin{defn}
A PMT $F$ is $\Bc$-\emph{structurally stable} if there exists a
holomorphic motion of $\Bc(F)$, parameterized by elements of a neighborhood
$\Nc_{F}\subset X_{PCM,K}$.
\end{defn}
As expected, the analogous result for rational maps also holds true
for PMTs.
\begin{thm}
Let $F$ be a structurally stable PMT; then it is $\Bc$-structurally
stable.
\end{thm}
\begin{proof}
Suppose that $F$ is not $\Bc$-structurally stable. Then, given a
holomorphic family $F_{\mu,\lambda}:\Chat\,\selfarrow$ parametrized
on $\Nc_{F}\subset X_{PCM,K}$, there does not exist a holomorphic
motion $\varphi_{\mu,\lambda}:\Bc(F)\rightarrow\Chat$ such that $\varphi_{\mu,\lambda}$
respects the dynamics in $\Bc(F)-B(F)$, or $\varphi_{\mu,\lambda}(B(F))\neq B(F_{\mu,\lambda})$,
for parameters close to $F$. In any case, $F$ and $F_{\mu,\lambda}$
cannot be topologically conjugated, and thus $F$ is not structurally
stable.
\end{proof}

\section{Structural Stability}

For rational maps, hyperbolic (or equivalently expanding) maps are
structurally stable. However, for PMTs, this is not the case, as reviewed
in Section \ref{sec:Hyperbolicity-and-Expansivity}.

On the other hand, we have the following
\begin{conjecture}
\label{thrm:Conj1}Let $F$ be a structurally stable PMT; then it
is hyperbolic and $\alpha$-expanding.
\end{conjecture}
\begin{rem}
Clearly, a structurally stable PMT cannot have parabolic, elliptic,
or identity periodic points, nor ghost-periodic points, because under
perturbations can be converted to attracting or repelling points.
The challenges in proving the previous conjecture arise from the following
cases of PMTs: i) those without periodic points where every regular
component is wandering, ii) the case with the pre-discontinuity set
dense in the sphere, or iii) the case with wandering components and
the pre-discontinuity set dense in some region with positive area.
\end{rem}
In line with the previous conjecture, the following can be proven:
\begin{thm}
Let $F$ be a structurally stable PMT without wandering domains; then
it is hyperbolic.
\end{thm}
\begin{proof}
Assume that $F$ is not hyperbolic and without wandering domains.
Then at least one of the following occurs:
\begin{enumerate}
\item $F$ has a parabolic, elliptic, or identity periodic point $z$. Under
perturbation of the component functions $f_{k}$ of $F$, $z$ can
be converted into an attracting or repelling periodic point for the
corresponding perturbed PMT $F_{\varepsilon}$.
\item $F$ has a ghost-periodic point $z$. Under perturbation of the discontinuity
set $B$, $z$ can be converted into a periodic point for the corresponding
perturbed PMT $F_{\varepsilon}$.
\item $\Bc(F)$ contains a region $U$ of positive area and $\Per(F)=\emptyset$.
\begin{enumerate}
\item If there exists a point $z\in\partial R_{i}\cap\partial R_{j}\cap U\subset B\cap U$,
then for every neighborhood $\Nc_{z}\subset U$ exists $w\in F^{-M}(B)\cap\Nc_{z}$
for some $M>0$, because of the density of $\left(\bigcup_{N\geq0}F^{-N}(B)\right)\cap U$
in $U$. Additionally, we can assume $w\in F^{-M}(B)\cap\Nc_{z}\subset R_{j}$.
Then a perturbation of $B$ around $F^{M}(w)$ (and possibly also
a perturbation of the component functions $f_{i}$ and $f_{j}$),
can lead to $F_{\varepsilon}^{-M}(B_{\varepsilon})\cap\Nc_{z}\cap R_{i}\neq\emptyset$,
where $F_{\varepsilon}$ is the corresponding perturbed PMT with $B(F_{\varepsilon})=B_{\varepsilon}$.
\item If there exists a point $z\in F^{-N}(B)\cap U$ with $N>0$ and $z\in R_{k}$
for some $k$, then for every neighborhood $\Nc_{z}\subset U\cap R_{k}$
exists $w\in F^{-M}(B)\cap\Nc_{z}$, for some $M>0$. Let $L=\min\left\{ N,M\right\} $,
$z_{0}=F^{L}(z)$ and $w_{0}=F^{L}(w)$. Then, $z_{0}\in B$ or $w_{0}\in B$,
and they are close to each other. Hence, we have sub-case (a).
\end{enumerate}
\end{enumerate}
In each of the three cases, $F$ cannot be topologically conjugated
with its corresponding perturbed $F_{\varepsilon}$.

The “without wandering domains” hypothesis ensures that the only instance
of $F$ such that $\Per(F)=\emptyset$ is the case (3) from the previous
list. 
\end{proof}
To ensure the structural stability of a PMT, several conditions must
be satisfied.
\begin{thm}
\label{thrm:GenStrucStab}Let $F$ be a PMT. If
\begin{enumerate}
\item each component function $f_{k}$ is loxodromic,
\item $F$ is hyperbolic and $\alpha$-expanding, and
\item $F$ is $\Bc$-structurally stable,
\end{enumerate}
then $F$ is structurally stable.
\end{thm}
\begin{proof}
By hypothesis (3), there exists a holomorphic family $F_{\mu,\lambda}:\Chat\,\selfarrow$
parameterized on $\Nc_{F}\subset X_{PCM,K}$, and a holomorphic motion
$\varphi_{\mu,\lambda}:\Bc(F)\rightarrow\Chat$ such that $\varphi_{\mu,\lambda}$
respects the dynamics in $\Bc(F)-B(F)$ and $\varphi_{\mu,\lambda}(B(F))=B(F_{\mu,\lambda})$.

Due to hypotheses (1) and (2), a possibly smaller neighborhood $\Nc_{F}$
can be chosen such that each $G\in\Nc_{F}$ satisfies hypotheses (1)
and (2), meaning $\alpha(G)$ contains all repelling but no parabolic
periodic points, and $\mathcal{R}(G)$ contains all attracting but
neither elliptic nor identity periodic points. Note that such PMTs
$G$ are constructed with the discontinuity set $B(G)=\varphi_{\mu,\lambda}(B(F))$,
and the component transformations $(g_{1},\dots,g_{K})$ are determined
by $(\mu,\lambda)\in\mathcal{N}_{F}$.

Utilizing the Bers-Royden extension, $\varphi$ has an extension to
a holomorphic motion $\Phi$ of $\Chat$ such that for each $(\mu,\lambda)\in\mathcal{N}_{F}$,
the function $h_{\mu,\lambda}=\Phi(\mu,\lambda,\_)$ is the unique
quasi-conformal homeomorphism on $\Chat$ with a harmonic Beltrami
differential in $\Chat-\Bc(F)$. Refer to the proof of Theorem \ref{thrm:StrucStabPSL2CK}
for further details regarding the construction of this extension.

By construction, $h_{\mu,\lambda}$ conjugates $F$ with $G$:
\begin{itemize}
\item If $z\in\Bc(F)-B$, by the definition of holomorphic motion of $\Bc(F)$,
we have $h_{\mu,\lambda}\circ F(z)=G\circ h_{\mu,\lambda}(z)$.
\item The function in $\mathcal{N}_{F}\times\Chat$ given by
\[
\widetilde{h}_{\lambda}(z)=\begin{cases}
g_{k}^{-1}\circ h_{\lambda}\circ f_{k} & \mathrm{if}\,z\in R_{k}\\
z & \mathrm{if}\,z\in B
\end{cases}
\]
is also an extension of the holomorphic motion $\varphi$, with a
harmonic Beltrami differential since $f_{k}$ and $g_{k}^{-1}$ are
holomorphic. By the uniqueness of the Bers-Royden extension under
such conditions, we have $\widetilde{h}_{\lambda}=h_{\lambda}$. 

Therefore, if $z\in\Rc(F)$, then $z\in R_{k}$ for some $k$, and
we can conclude
\[
h_{\lambda}\circ F(z)=h_{\lambda}\circ f_{k}(z)=g_{k}\circ h_{\lambda}(z)=G\circ h_{\lambda}(z).
\]

\end{itemize}
\end{proof}
Based on experimental evidence, the equivalence between structural
stability and the conditions of the previous theorem appears to be
true. Hence, a stronger conjecture than conjecture \ref{thrm:Conj1}
above is: 
\begin{conjecture}
If $F$ is a structurally stable PMT, then each component transformation
$f_{k}$ is loxodromic, $F$ is hyperbolic, and $F$ is $\alpha$-expanding.
\end{conjecture}

\section{Example: The Tent Maps Family}

To finalize the analysis of the stability of PMTs, we will demonstrate
applications of previous results to the complex version of the well-known
family of tent maps in $\R$.
\begin{defn}
The family of \emph{complex tent maps} 
\[
\left\{ T_{B,\lambda}:\Chat\,\selfarrow\right\} _{B\in\Pc_{2},\,\lambda\in\C-\left\{ 0\right\} }
\]
is defined by
\[
T_{B,\lambda}(z)=\begin{cases}
f_{1}(z) & \mathrm{if}\,z\in R_{1}\\
f_{2}(z) & \mathrm{if}\,z\in R_{2},
\end{cases}
\]
where $f_{1}(z)=\lambda z$, $f_{2}(z)=\lambda-\lambda z$, $B=\partial R_{1}=\partial R_{2}$
and $\frac{1}{2}\in B$.
\end{defn}
\begin{rem}
The condition $\frac{1}{2}\in B$ s necessary to ensure similar behavior
to the real case: $f_{1}(\frac{1}{2})=f_{2}(\frac{1}{2})=\lambda\frac{1}{2}$.
However, $T_{B,\lambda}$ cannot be extended to a continuous function
in every neighborhood $\Nc_{\frac{1}{2}}$.
\end{rem}
Let us list several facts about this family of maps.
\begin{itemize}
\item Clearly, it is a holomorphic family of PMTs.
\item The fixed points of $f_{1}$ are $0$ and $\infty$. The fixed points
of $f_{2}$ are $z_{\lambda}=\frac{\lambda}{\lambda+1}$ and $\infty$.
Thus
\[
\Fix(T_{B,\lambda})=\left(\left(\left\{ 0,\infty\right\} \cap R_{1}\right)\cup\left(\left\{ z_{\lambda},\infty\right\} \cap R_{2}\right)\right)\cap\left(\Rc(T_{B,\lambda})\cup\alpha(T_{B,\lambda})\right).
\]
\item If $|\lambda|<1$, then $f_{1}$ and $f_{2}$ are affine contractions
in $\C$. Therefore, for almost every $\lambda\in\D$, all points
in $\Rc(F)$ tend to an attracting or to a ghost periodic orbit. Also,
it can be shown that if $B\subset\C$, then $\alpha(T_{B,\lambda})=\left\{ \infty\right\} $
(see \cite{Ler2016}).
\item If $|\lambda|=1$, then $f_{1}$ and $f_{2}$ are euclidean isometries.
If $B\subset\C$, then every point in $\Rc(F)$ is periodic or pre-periodic
(see \cite{Goe2000} for this result).
\item If $\lambda=1$, then $f_{1}=Id|_{R_{1}}$ and $f_{2}$ is a euclidean
rotation. If $\lambda=-1$, then $f_{1}$ is a euclidean rotation
and $f_{2}$ is a translation. In any case, every point in $\Rc(F)$
is periodic or pre-periodic (see \cite{Ler2016}).
\item If $|\lambda|>1$ and $B\subset\C$, then $\infty$ is an attracting
fixed point of $T_{B,\lambda}$.
\end{itemize}
The global behavior of the orbits can be determined when parameters
are such that $|\lambda|\neq1$ (see \cite{Ler2016}).
\begin{thm}
~
\begin{itemize}
\item If $|\lambda|<1$, $T_{B,\lambda}$ is globally attracting; that is,
there exists $r\in(0,\infty)$ such that if $z\in\Rc(T_{B,\lambda})-\left\{ \infty\right\} $,
then there exists $N\in\N$ such that \textbar$T_{B,\lambda}^{n}(z)|\leq r$
for all $n\geq N$.
\item If $|\lambda|>1$, $T_{B,\lambda}$ is globally repelling; that is,
there exists $r\in(0,\infty)$ such that if $|z|>r$ and $z\in\Rc(T_{B,\lambda})$,
then $\underset{n\rightarrow\infty}{\lim}T_{B,\lambda}^{n}(z)=\infty$.
\end{itemize}
\end{thm}
\begin{figure}[h]
\begin{centering}
\includegraphics{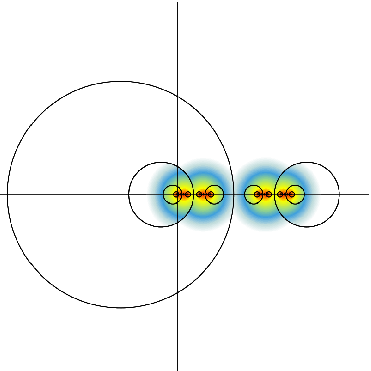}\hspace*{0.25cm}\includegraphics{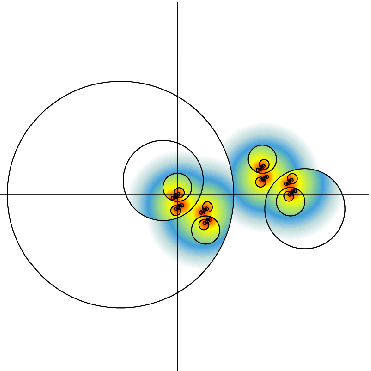}
\par\end{centering}
\begin{centering}
\includegraphics{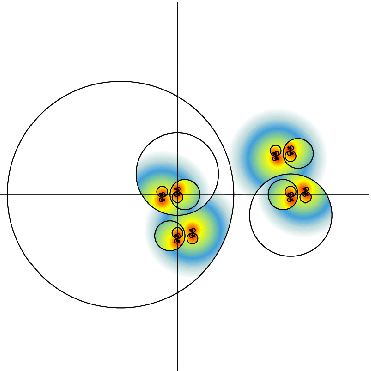}\hspace*{0.25cm}\includegraphics{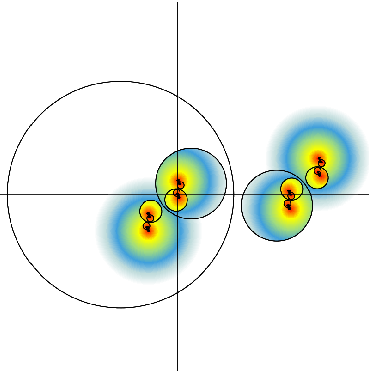}
\par\end{centering}
\caption{\protect\label{fig:Fig5}Pre-discontinuity and regular sets of the
tent maps $T_{B,\lambda}$ from Example \ref{example:TentSS}.\protect \\
Top left: With $\lambda=\frac{7}{2}$. Top right: With $\lambda=2+2i$.\protect \\
Bottom left: With $\lambda=\frac{11}{4}i$. Bottom right: With $\lambda=-\frac{5}{2}+2i$.}
\end{figure}

Notice that for parameters such that $|\lambda|\neq1$, $f_{1}$ and
$f_{2}$ are loxodromic and $\Fix(f_{1})\cap\Fix(f_{2})=\left\{ \infty\right\} $,
then the group $\Gamma=\big<f_{1},f_{2}\big>$ is not discrete. Similarly,
when $\lambda=e^{2\pi\theta i}$ with $\theta$ being an irrational
number, $\Gamma=\big<f_{1},f_{2}\big>$ is not discrete. In any case,
we have that $\Lambda(\Gamma)$ (the limit set of $\Gamma$) is $\Chat$,
and the results regarding stability related to structurally stable
Kleinian groups cannot be applied (see \cite{Ler2005,LerSie2019,Ler2023}
for these results).

However, structural stability can be found in the family under the
following conditions:
\begin{enumerate}
\item Parameter $|\lambda|\neq1$.
\item Bounded discontinuity set, that is, $B\subset\C$.
\item Finite fixed points ($0$ and $z_{\lambda}$) of $f_{1}$ and $f_{2}$
such that they are not in $B$.
\item Pre-discontinuity set formed exclusively by homeomorphic copies of
$B$ and the corresponding $\alpha$-limit set. This can be achieved
by selecting $\lambda$ with a sufficiently large or a sufficiently
small modulus.
\end{enumerate}
Then, we have
\begin{itemize}
\item By (1), $f_{1}$ and $f_{2}$ are loxodromic.
\item $T_{B,\lambda}$ has no ghost-fixed points, since $\infty,0,z_{\lambda}\notin B$
by incises (2) and (3).
\item $\infty$ is an attracting or repelling fixed point of $T_{B,\lambda}$,
according to (1) and (2).
\item $T_{B,\lambda}$ is hyperbolic and $\alpha$-expanding. By (1) and
(4):
\begin{itemize}
\item Every point in $\alpha(T_{B,\lambda})$ is either a repelling periodic
point, a pre-repelling periodic point, or has an infinite orbit while
being a limit point of the semi-group generated by $f_{1}^{-1}$ and
$f_{2}^{-1}$.
\item Every point in $\Rc(T_{B,\lambda})$ is attracted to $\infty$ when
$|\lambda|>1$, or to $0$ (if $0\in R_{1}$), or to $z_{\lambda}$
(if $z_{\lambda}\in R_{2}$), when $|\lambda|<1$. 
\end{itemize}
\end{itemize}
In summary, $T_{B,\lambda}$ fulfilling (1), (2), (3), and (4) has
loxodromic component transformations, is hyperbolic, and is $\alpha$-expanding.
Clearly, a holomorphic motion can be constructed for each $\Bc(T_{B,\lambda})$,
and then, by Theorem \ref{thrm:GenStrucStab}, all these PMTs $T_{B,\lambda}$
are structurally stable.
\begin{example}
\label{example:TentSS}The pre-discontinuity sets of $T_{B,\lambda}$
with $R_{1}=\left\{ z:\,|z+\frac{1}{2}|<1\right\} $ and $R_{2}=\Chat-\overline{R_{1}}$
are depicted in black in the images from Figure \ref{fig:Fig5}. The
color gradient indicates the proximity of repelling periodic points
in $\alpha(T_{B,\lambda})$.

In this example, $B$ is fixed; however, it is clear that such $B$
can be deformed while preserving the conditions (1), (2), (3), and
(4) mentioned above. Consequently, these new maps generated by deformation
are structurally stable.

\rule[1ex]{0.1\paperwidth}{0.25pt}
\end{example}

\end{document}